\documentclass{birkmult}

\usepackage[sort]{cite}       % nicer citations
\usepackage{amssymb}          % various ams stuff
\usepackage{latexsym}         % dito
\usepackage{amsmath}          % dito
\usepackage[latin1]{inputenc} % smart input of special characters
\usepackage{eucal}            % nicer caligraphic fonts
\usepackage{bbm}              % nicer bbm fonts
\usepackage{mathrsfs}         % nette kaligraphische Schrift
\usepackage{stmaryrd}         % noch ein paar mehr Symbole (Klammern)

\renewcommand{\mathbb}[1]{\mathbbm{#1}} % use nicer bbm fonts

\newcommand{\Lie}        {\operatorname{\mathscr{L}}}    
\newcommand{\cc}[1]      {\overline{{#1}}}              
\newcommand{\id}         {\operatorname{\mathsf{id}}}   
\newcommand{\Hom}        {\operatorname{\mathsf{Hom}}}   
\newcommand{\End}        {\operatorname{\mathsf{End}}}   
\newcommand{\SP}[1]      {\left\langle{#1}\right\rangle} 
\newcommand{\ver}        {\mathrm{ver}}                
\newcommand{\pr}         {\mathrm{pr}}                    
\newcommand{\I}          {\mathrm{i}}
\newcommand{\E}          {\mathrm{e}}
\newcommand{\D}          {\operatorname{\mathrm{d}}} 
\newcommand{\Anti}       {\Lambda}
\newcommand{\Schouten}[1]{\left\llbracket{#1}\right\rrbracket}
\newcommand{\Diffop}     {\operatorname{\mathrm{Diffop}}}
\newcommand{\ring}[1]    {\mathsf{#1}}
\newcommand{\HCdiff}     {\operatorname{\mathrm{HC}}_{\mathrm{diff}}}
\newcommand{\HHdiff}     {\operatorname{\mathrm{HH}}_{\mathrm{diff}}}
\newcommand{\Def}        {\operatorname{\mathrm{Def}}}

\newtheorem{theorem}{Theorem}[section]
\newtheorem{corollary}[theorem]{Corollary}

\theoremstyle{definition}
\newtheorem{definition}[theorem]{Definition}
\theoremstyle{remark}
\newtheorem{remark}[theorem]{Remark}
\newtheorem*{example}{Example}

\numberwithin{equation}{section}

\begin{document}

%-------------------------------------------------------------------------
% editorial commands: to be inserted by the editorial office
%
%\firstpage{1}
%\volume{228}
%\Copyrightyear{2004}
%\DOI{003-0001}
%
%
%\seriesextra{Just an add-on}
%\seriesextraline{This is the Concrete Title of this Book\br H.E. R and S.T.C. W, Eds.}
%
% for journals:
%
%\firstpage{1}
%\issuenumber{1}
%\Volumeandyear{1 (2004)}
%\Copyrightyear{2004}
%\DOI{003-xxxx-y}
%\Signet
%\commby{inhouse}
%\submitted{March 14, 2003}
%\received{March 16, 2000}
%\revised{June 1, 2000}
%\accepted{July 22, 2000}
%
%
%
%---------------------------------------------------------------------------

%
% Title author etc...
%

\title{Noncommutative field theories from a deformation point of view}

\author{Stefan Waldmann}

\address{%
  Fakult{\"a}t f{\"u}r Mathematik und Physik \\
  Albert-Ludwigs-Universit{\"a}t Freiburg \\
  Physikalisches Institut \\
  Hermann Herder Strasse 3 \\
  D 79104 Freiburg \\
  Germany
}

\email{Stefan.Waldmann@physik.uni-freiburg.de}

%
% classification, keywords, date
%

\subjclass{Primary 53D55; Secondary 58B34, 81T75}

\keywords{Noncommutative field theory, Deformation quantization,
  Principal Bundles}

\date{October 2007}

\begin{abstract}
    In this review we discuss the global geometry of noncommutative
    field theories from a deformation point of view: The space-times
    under consideration are deformations of classical space-time
    manifolds using star products. Then matter fields are encoded in
    deformation quantizations of vector bundles over the classical
    space-time. For gauge theories we establish a notion of
    deformation quantization of a principal fiber bundle and show how
    the deformation of associated vector bundles can be obtained.
\end{abstract}

%
% Title
%
\maketitle

%
% Introduction
%

\section{Introduction}
\label{sec:Introduction}

Noncommutative geometry is commonly believed to be a reasonable
candidate for the marriage of classical gravity theory in form of
Einstein's general relativity on one hand and quantum theory on the
other hand. Both theories are experimentally well-established within
large regimes of energy and distance scales. However, from a more
fundamental point of view, the coexistence of these two theories
becomes inevitably inconsistent when one approaches the Planck scale
where gravity itself gives significant quantum effects.

Since general relativity is ultimately the theory of the geometry of
space-time it seems reasonable to use notions of `quantum geometry'
known under the term \emph{noncommutative geometry} in the sense of
Connes \cite{connes:1994a} to achieve appropriate formulations of what
eventually should become quantum gravity. Of course, this ultimate goal
has not yet been reached but techniques of noncommutative geometry
have been used successfully to develop models of quantum field
theories on quantum space-times being of interest for their own.
Moreover, a deeper understanding of ordinary quantum field theories
can be obtained by studying their counterparts on `nearby'
noncommutative space-times. On the other hand, people started to
investigate experimental implications of a possible  noncommutativity
of space-time in future particle experiments.

Such a wide scale of applications and interests justifies a more
\emph{conceptual} discussion of noncommutative space-times and
(quantum) field theories on them in order to clarify fundamental
questions and generic features expected to be common to all examples.

In this review, we shall present such an approach from the point of
view of deformation theory: noncommutative space-times are not studied
by themselves but always with respect to a classical space-time, being
suitably deformed into the noncommutative one. Clearly, this point of
view can not cover all possible (and possibly interesting)
noncommutative geometries but only a particular class. Moreover, we
focus on \emph{formal} deformations for technical reasons. It is
simply the most easy approach where one can rely on the very powerful
machinery of algebraic deformation theory. But it also gives hints on
approaches beyond formal deformations: finding obstructions in the
formal framework will indicate even more severe obstructions in any
non-perturbative approach.

In the following, we discuss mainly two questions: first, what is the
appropriate description of matter fields on deformed space-times and,
second, what are the deformed analogues of principal bundles needed
for the formulation of gauge theories. The motivation for these two
questions should be clear.

The review is organized as follows: in Section~\ref{sec:NCSpacetimes},
we recall some basic definitions and properties concerning deformation
quantizations and star products needed for the set-up of
noncommutative space-times. We discuss some fundamental examples as
well as a new class of locally noncommutative space-times.
Section~\ref{sec:MatterFieldsVectorBundles} is devoted to the study of
matter fields: we use the Serre-Swan theorem to relate matter fields
to projective modules and discuss their deformation theory. Particular
interest is put on the mass terms and their positivity properties.  In
Section~\ref{sec:DeformedPrincipalBundles} we establish the notion of
deformation quantization of principal fiber bundles and discuss the
existence and uniqueness results. Finally, in
Section~\ref{sec:CommutantAssociated} we investigate the resulting
commutant and formulate an appropriate notion of associated (vector)
bundles. This way we make contact to the results of
Section~\ref{sec:MatterFieldsVectorBundles}. The review is based on
joint works with Henrique Bursztyn on one hand as well as with Martin
Bordemann, Nikolai Neumaier and Stefan Weiß on the other hand.

%
% Noncommutative space-times
%

\section{Noncommutative space-times}
\label{sec:NCSpacetimes}

In order to implement uncertainty relations for measuring coordinates
of events in space-time it has been proposed already very early to
replace the commutative algebra of (coordinate) functions by some
noncommutative algebra. In \cite{doplicher.fredenhagen.roberts:1995a}
a concrete model for a noncommutative Minkowski space-time was
introduced with commutation relations of the form
\begin{equation}
    \label{eq:xmuxnuthetamunu}
    [\hat{x}^\mu, \hat{x}^\nu] = \I \lambda \theta^{\mu\nu},
\end{equation}
where $\lambda$ plays the role of the deformation parameter and has
the physical dimension of an area. Usually, this area will be
interpreted as the Planck area. Moreover, $\theta$ is a real,
antisymmetric tensor which in
\cite{doplicher.fredenhagen.roberts:1995a} and many following papers
is assumed to be \emph{constant}: in
\cite{doplicher.fredenhagen.roberts:1995a} this amounts to require
that $\theta^{\mu\nu}$ belongs to the center of the new algebra of
noncommutative coordinates.

Instead of constructing an abstract algebra where commutation
relations like \eqref{eq:xmuxnuthetamunu} are fulfilled, it is
convenient to use a `symbol calculus' and encode
\eqref{eq:xmuxnuthetamunu} already for the classical coordinate
functions by changing the multiplication law instead. For functions
$f$ and $g$ on the classical Minkowski space-time one defines the
Weyl-Moyal star product by
\begin{equation}
    \label{eq:WeylMoyal}
    f \star g = \mu \circ \E^{\frac{\I\lambda}{2} \theta^{\mu\nu}
      \frac{\partial}{\partial x^\mu} \otimes
      \frac{\partial}{\partial x^\nu}}
    (f \otimes g),
\end{equation}
where $\mu (f \otimes g) = fg$ denotes the undeformed, pointwise
product. Then \eqref{eq:xmuxnuthetamunu} holds for the classical
coordinate functions with respect to the $\star$-commutator.

Clearly, one has to be slightly more careful with expressions like
\eqref{eq:WeylMoyal}: in order to make sense out of the infinite
differentiations the functions $f$ and $g$ first should be $C^\infty$.
But then the exponential series does not converge in general whence a
more sophisticated analysis is required. Though this can be done in a
completely satisfying way for this particular example, we shall not
enter this discussion here but consider \eqref{eq:WeylMoyal} as a
formal power series in the deformation parameter $\lambda$. Then
$\star$ becomes an associative $\mathbb{C}[[\lambda]]$-bilinear
product for $C^\infty(\mathbb{R}^4)[[\lambda]]$, i. e. a star product
in the sense of \cite{bayen.et.al:1978a}. It should be noted that the
interpretation of \eqref{eq:WeylMoyal} as formal series in $\lambda$
is physically problematic: $\lambda$ is the Planck area and hence a
physically measurable and non-zero quantity. Thus our point of view
only postpones the convergence problem and can be seen as a
perturbative approach.

With this example in mind, one arrives at several conceptual
questions: The first is that Minkowski space-time is clearly not a
very realistic background when one wants to consider quantum effects
of `hard' gravity. Here already classically nontrivial curvature and
even nontrivial topology may arise. Thus one is forced to consider
more general and probably even generic Lorentz manifolds instead.
Fortunately, deformation quantization provides a well-established and
successful mathematical framework for this geometric situation.

Recall that a star product on a manifold $M$ is an associative
$\mathbb{C}[[\lambda]]$-bilinear multiplication $\star$ for $f, g \in
C^\infty(M)[[\lambda]]$ of the form
\begin{equation}
    \label{eq:StarProduct}
    f \star g = \sum_{r=0}^\infty \lambda^r C_r(f, g),
\end{equation}
where $C_0(f, g) = fg$ is the undeformed, pointwise multiplication and
the $C_r$ are bidifferential operators. Usually, one requires $1 \star
f = f = f \star 1$ for all $f$. It is easy to see that $\{f, g\} =
\frac{1}{\I} (C_1(f, g) - C_1(g, f))$ defines a Poisson bracket on
$M$. Conversely, and this is the highly nontrivial part, any Poisson
bracket $\{f, g\} = \theta(\D f, \D g)$, where
\begin{equation}
    \label{eq:ThetaGlobal}
    \theta \in \Gamma^\infty(\Anti^2 TM),
    \quad
    \Schouten{\theta, \theta} = 0
\end{equation}
is the corresponding Poisson tensor, can be quantized into a star
product \cite{kontsevich:2003a, dewilde.lecomte:1983b}. Beside these
existence results one has a very good understanding of the
classification of such star products \cite{nest.tsygan:1995a,
  gutt.rawnsley:1999a, kontsevich:2003a}, see also
\cite{dito.sternheimer:2002a, gutt:2000a} for recent reviews and
\cite{waldmann:2007a} for an introduction.

With this geometric interpretation the Weyl-Moyal star product on
Minkowski space-time turns out to be a deformation quantization of the
\emph{constant} Poisson structure
\begin{equation}
    \label{eq:ThetaConstant}
    \theta = \frac{1}{2} \theta^{\mu\nu} 
    \frac{\partial}{\partial x^\mu}
    \wedge
    \frac{\partial}{\partial x^\nu}.
\end{equation}
On a generic space-time $M$ there is typically \emph{no} transitive
action of isometries which would justify the notion of a `constant'
bivector field. Thus a star product $\star$ on $M$ is much more
complicated than \eqref{eq:WeylMoyal} in general: already the first
order term is a (nontrivial) Poisson structure and for the higher
order terms one has to invoke the (unfortunately rather inexplicit)
existence theorems.

Thus answering the first question by using general star products
raises the second: what is the physical role of a Poisson structure on
space-time? While on Minkowski space-time with constant $\theta$ we
can view the finite number of coefficients $\theta^{\mu\nu} \in
\mathbb{R}$ as \emph{parameters} of the theory this is certainly no
longer reasonable in the more realistic geometric framework: there is
an infinity of Poisson structures on each manifold whence an
interpretation as `parameter' yields a meaningless theory. Instead,
$\theta$ has to be considered as a \emph{field} itself, obeying its
own dynamics compatible with the constraint of the Jacobi identity
$\Schouten{\theta, \theta} = 0$. Unfortunately, up to now a reasonable
`field equation' justified by first principles seems to be missing.

This raises a third conceptual question, namely why should there by
any Poisson structure on $M$ and what are possible experimental
implications? In particular, the original idea of introducing a
noncommutative structure was to implement uncertainty relations
forbidding the precise localization of events. The common believe is
that such quantum effects should only play a role when approaching the
Planck scale. Now it turns out that the quantum field theories put on
such a noncommutative Minkowski space-time (or their Euklidian
counterparts) suffer all from quite unphysical properties: Typically,
the noncommutativity enters in long-distance/low-energy features
contradicting our daily life experience. Certainly, a last word is not
said but there might be a simple explanation why such effects should
be expected: the \emph{global} $\theta$ (constant or not) yields
global effects on $M$. This was the starting point of a more refined
notion of noncommutative space-times advocated in
\cite{bahns.waldmann:2007a, heller.neumaier.waldmann:2007a} as
\emph{locally} noncommutative space-times. Roughly speaking, without
entering the technical details, it is not $M$ which should become
noncommutative but $TM$. Here the tangent bundle is interpreted as the
bundle of all normal charts on $M$ and for each normal chart with
origin $p \in M$ one constructs its own star product $\star_p$. The
crucial property is then that $\star_p$ is the pointwise, commutative
product outside a (small) compact subset around $p$. This way, the
long-distance behaviour (with respect to the reference point $p$) is
classical while close to $p$ there is a possibly even very strong
noncommutativity. In some sense, this is an implementation of an idea
of Julius Wess, proposing that the transition from classical geometry
to quantum geometry should be understood as a kind of phase transition
taking place at very small distances \cite{wess:misc}. Of course, the
conceptual question about the physical origin of the corresponding
Poisson structure on $TM$ as well as the convergence problem still
persists also in this approach.

Ignoring these questions about the nature of $\theta$, we shall assume
in the following that we are given a star product $\star$ on a
manifold $M$ which can be either space-time itself or its tangent
bundle in the locally noncommutative case. Then we address the
question how to formulate reasonable field theories on $(M, \star)$.
Here we shall focus on \emph{classical} field theories which still
need to be quantized later on. On the other hand, we seek for a
\emph{geometric} formulation not relying on particular assumptions
about the underlying classical space-time.

%
% Matter fields and deformed vector bundles
%

\section{Matter fields and deformed vector bundles}
\label{sec:MatterFieldsVectorBundles}

In this section we review some results from
\cite{bursztyn.waldmann:2005b, waldmann:2005a, waldmann:2001b,
  bursztyn.waldmann:2000b}.

In classical field theories both bosonic and fermionic matter fields
are given by sections of appropriate vector bundles. For convenience,
we choose the vector bundles to be complex as also the function
algebra $C^\infty(M)$ consists of complex-valued functions. However,
the real case can be treated completely analogously.  Thus let $E
\longrightarrow M$ be a complex vector bundle over $M$. Then the
$E$-valued fields are the (smooth) sections $\Gamma^\infty(E)$ which
form a module over $C^\infty(M)$ by pointwise multiplication. Thanks
to the commutativity of $C^\infty(M)$ we have the freedom to choose
this module structure to be a right module structure for later
convenience.

It is a crucial feature of vector bundles that $\Gamma^\infty(E)$ is
actually a finitely generated and projective module:
\begin{theorem}[Serre-Swan]
    \label{theorem:SerreSwan}
    The sections $\Gamma^\infty(E)$ of a vector bundle $E
    \longrightarrow M$ are a finitely generated and projective
    $C^\infty(M)$-module. Conversely, any such module arises this way
    up to isomorphism.
\end{theorem}
Recall that a right module $\mathcal{E}_{\mathcal{A}}$ over an algebra
$\mathcal{A}$ is called finitely generated and projective if there
exists an idempotent $e^2 = e \in M_n(\mathcal{A})$ such that
$\mathcal{E}_{\mathcal{A}} \cong e\mathcal{A}^n$ as right
$\mathcal{A}$-modules. More geometrically speaking, for any vector
bundle $E \longrightarrow M$ there is another vector bundle $F
\longrightarrow M$ such that their Whitney sum $E \oplus F$ is
isomorphic to a trivial vector bundle $M \times \mathbb{C}^n
\longrightarrow M$. Note that the Serre-Swan theorem has many
incarnations, e.g. the original version was formulated for compact
Hausdorff spaces and continuous sections/functions. Note also that for
our situation no compactness assumption is necessary (though it
drastically simplifies the proof) as manifolds are assumed to be
second countable.
\begin{remark}
    \label{remark:SerreSwan}
    The Serre-Swan theorem is the main motivation for noncommutative
    geometry to consider finitely generated and projective modules
    over a not necessarily commutative algebra $\mathcal{A}$ as
    `vector bundles' over the (noncommutative) space described by
    $\mathcal{A}$ in general.
\end{remark}

For physical applications in field theory one usually has more
structure on $E$ than just a bare vector bundle. In particular, for a
Lagrangean formulation a `mass term' in the Lagrangean is
needed. Geometrically such a mass term corresponds to a Hermitian
fiber metric $h$ on $E$. One can view a Hermitian fiber metric as a
map
\begin{equation}
    \label{eq:FiberMetric}
    h: \Gamma^\infty(E) \times \Gamma^\infty(E) 
    \longrightarrow C^\infty(M),
\end{equation}
which is $\mathbb{C}$-linear in the second argument and satisfies
$\cc{h(\phi, \psi)} = h(\psi, \phi)$, $h(\phi, \psi f) = h(\phi, \psi)
f$ as well as 
\begin{equation}
    \label{eq:Positive}
    h(\phi, \phi) \ge 0    
\end{equation}
for $\phi, \psi \in \Gamma^\infty(E)$ and $f \in C^\infty(M)$. The
pointwise non-degeneracy of $h$ is equivalent to the property that
\begin{equation}
    \label{eq:hNondegenerate}
    \Gamma^\infty(E) \ni \phi 
    \; \mapsto \; 
    h(\phi, \cdot) \in \Gamma^\infty(E^*)
\end{equation}
is an antilinear module isomorphism. Note that the sections of the
dual vector bundle $E^* \longrightarrow M$ coincide with the dual
module, i.e. we have $\Gamma^\infty(E^*) = \Hom_{C^\infty(M)}
(\Gamma^\infty(E), C^\infty(M))$.

In order to encode now the positivity \eqref{eq:Positive} in a more
algebraic way suitable for deformation theory, we have to consider the
following class of algebras: First, we use a ring of the form
$\ring{C} = \ring{R}(\I)$ with $\I^2 = -1$ for the scalars where
$\ring{R}$ is an \emph{ordered ring}. This includes both $\mathbb{R}$
and $\mathbb{R}[[\lambda]]$, where positive elements in
$\mathbb{R}[[\lambda]]$ are defined by
\begin{equation}
    \label{eq:PositivePowerseries}
    a = \sum_{r=r_0}^\infty \lambda^r a_r > 0
    \quad
    \textrm{if}
    \quad a_{r_0} > 0.
\end{equation}
In fact, this way $\ring{R}[[\lambda]]$ becomes an ordered ring
whenever $\ring{R}$ is ordered. More physically speaking, the ordering
of $\mathbb{R}[[\lambda]]$ refers to a kind of asymptotic
positivity. Then the algebras in question should be $^*$-algebras over
$\ring{C}$: Indeed, $C^\infty(M)$ is a $^*$-algebra over $\mathbb{C}$
where the $^*$-involution is the pointwise complex conjugation. For
the deformed algebras $(C^\infty(M)[[\lambda]], \star)$ we require
that the star product is Hermitian, i.e.
\begin{equation}
    \label{eq:HermitianStarProduct}
    \cc{f \star g} = \cc{g} \star \cc{f}
\end{equation}
for all $f, g \in C^\infty(M)[[\lambda]]$. For a real Poisson
structure $\theta$ this can be achieved by a suitable choice of
$\star$.

For such a $^*$-algebra we can now speak of positive functionals and
positive elements \cite{bursztyn.waldmann:2001a} by mimicking the
usual definitions from operator algebras, see e.g.
\cite{schmuedgen:1990a} for the case of (unbounded) operator algebras
and \cite{waldmann:2004a} for a detailed comparison.
\begin{definition}
    \label{definition:Positive}
    Let $\mathcal{A}$ be a $^*$-algebra over $\ring{C} =
    \ring{R}(\I)$. A $\ring{C}$-linear functional $\omega: \mathcal{A}
    \longrightarrow \ring{C}$ is called positive if $\omega(a^*a) \ge
    0$ for all $a \in \mathcal{A}$. An element $a \in \mathcal{A}$ is
    called positive if $\omega(a) \ge 0$ for all positive functionals
    $\omega$.
\end{definition}
We denote the convex cone of positive elements in $\mathcal{A}$ by
$\mathcal{A}^+$. It is an easy exercise to show that for $\mathcal{A}
= C^\infty(M)$ the positive functionals are the compactly supported
Borel measures and $\mathcal{A}^+$ consists of functions $f$ with
$f(x) \ge 0$ for all $x \in M$.

Using this notion of positive elements and motivated by
\cite{lance:1995a}, the algebraic formulation of a fiber metric is now
as follows \cite{bursztyn.waldmann:2005b, bursztyn.waldmann:2000b}:
\begin{definition}
    \label{definition:FiberMetric}
    Let $\mathcal{E}_{\mathcal{A}}$ be a right
    $\mathcal{A}$-module. Then an inner product $\SP{\cdot, \cdot}$ on
    $\mathcal{E}_{\mathcal{A}}$ is a map
    \begin{equation}
        \label{eq:InnerProduct}
        \SP{\cdot, \cdot}: 
        \mathcal{E}_{\mathcal{A}} \times \mathcal{E}_{\mathcal{A}}
        \longrightarrow \mathcal{A},
    \end{equation}
    which is $\ring{C}$-linear in the second argument and satisfies
    $\SP{x, y} = \SP{y, x}^*$, $\SP{x, y \cdot a} = \SP{x, y} a$, and
    $\SP{x, y} = 0$ for all $y$ implies $x = 0$. The inner product is
    called strongly non-degenerate if in addition
    \begin{equation}
        \label{eq:StronglyNondegenerate}
        \mathcal{E}_{\mathcal{A}} \ni x
        \; \mapsto \; \SP{x, \cdot} 
        \in \mathcal{E}^* =
        \Hom_{\mathcal{A}}(\mathcal{E}_{\mathcal{A}}, \mathcal{A})
    \end{equation}
    is bijective. It is called completely positive if for all $n \in
    \mathbb{N}$ and $x_1, \ldots, x_n \in \mathcal{E}_{\mathcal{A}}$
    one has $(\SP{x_i, x_j}) \in M_n(\mathcal{A})^+$.
\end{definition}
Clearly, a Hermitian fiber metric on a complex vector bundle endows
$\Gamma^\infty(E)$ with a completely positive, strongly non-degenerate
inner product in the sense of Definition~\ref{definition:FiberMetric}.

With the above definition in mind we can now formulate the following
deformation problem \cite{bursztyn.waldmann:2000b}:
\begin{definition}
    \label{definition:DeformationOfE}
    Let $\star$ be a Hermitian star product on $M$ and $E
    \longrightarrow M$ a complex vector bundle with fiber metric $h$.
    \begin{enumerate}
    \item A deformation quantization $\bullet$ of $E$ is a right
        module structure $\bullet$ for $\Gamma^\infty(E)[[\lambda]]$
        with respect to $\star$ of the form
        \begin{equation}
            \label{eq:RightModulephif}
            \phi \bullet f = \sum_{r=0}^\infty \lambda^r R_r(\phi, f)
        \end{equation}
        with bidifferential operators $R_r$ and $R_0(\phi, f) = \phi
        f$.
    \item For a given deformation quantization $\bullet$ of $E$ a
        deformation quantization of $h$ is a completely positive inner
        product $\boldsymbol{h}$ for $(\Gamma^\infty(E)[[\lambda]],
        \bullet)$ of the form
        \begin{equation}
            \label{eq:DeformedFiberMetric}
            \boldsymbol{h}(\phi, \psi) = \sum_{r=0}^\infty \lambda^r
            \boldsymbol{h}_r(\phi, \psi)
        \end{equation}
        with (sesquilinear) bidifferential operators
        $\boldsymbol{h}_r$ and $\boldsymbol{h}_0 = h$.
    \end{enumerate}
\end{definition}

In addition, we call two deformations $\bullet$ and $\tilde{\bullet}$
\emph{equivalent} if there exists a formal series of differential
operators
\begin{equation}
    \label{eq:EquivalenceTrafo}
    T = \id + \sum_{r=1}^\infty \lambda^r T_t:
    \Gamma^\infty(E)[[\lambda]] \longrightarrow
    \Gamma^\infty(E)[[\lambda]],
\end{equation}
such that
\begin{equation}
    \label{eq:Equivalencebullets}
    T(\phi \bullet f) = T(\phi) \tilde{\bullet} f.
\end{equation}
With other words, $T$ is a module isomorphism starting with the
identity in order $\lambda^0$ such that $T$ is not visible in the
classical/commutative limit. Conversely, starting with one deformation
$\bullet$ and a $T$ like in \eqref{eq:EquivalenceTrafo}, one obtains
another equivalent deformation $\tilde{\bullet}$ by defining
$\tilde{\bullet}$ via \eqref{eq:Equivalencebullets}. Similarly, we
define two deformations $\boldsymbol{h}$ and $\tilde{\boldsymbol{h}}$
to be \emph{isometric} if there exists a self-equivalence $U$ with
\begin{equation}
    \label{eq:Isometric}
    \boldsymbol{h}(\phi, \psi)
    = \tilde{\boldsymbol{h}}(U(\phi), U(\psi)).
\end{equation}
The relevance of the above notions for noncommutative field theories
should now be clear: for a classical matter field theory modelled on
$E \longrightarrow M$ we obtain the corresponding noncommutative field
theory by choosing a deformation $\bullet$ (if existing!) together
with a deformation $\boldsymbol{h}$ (if existing!) in order to write
down noncommutative Lagrangeans involving expressions like
$\mathcal{L}(\phi) = \boldsymbol{h}(\phi, \phi) + \cdots$.

Note that naive expressions like $\cc{\phi} \star \phi$ do not make
sense geometrically, even on the classical level: sections of a vector
bundle can not be `multiplied' without the extra structure of a fiber
metric $h$ unless the bundle is trivial \emph{and} trivialized. In
this particular case we can of course use the canonical fiber metric
coming from the canonical inner product on $\mathbb{C}^n$. We refer to
\cite{waldmann:2005a, waldmann:2001b} for a further discussion.

We can now state the main results of this section, see
\cite{bursztyn.waldmann:2005b, bursztyn.waldmann:2000b} for detailed
proofs:
\begin{theorem}
    \label{theorem:ExistenceDeformedE}
    For any star product $\star$ on $M$ and any vector bundle $E
    \longrightarrow M$ there exists a deformation quantization
    $\bullet$ with respect to $\star$ which is unique up to
    equivalence.
\end{theorem}
\begin{theorem}
    \label{theorem:ExistenceFiberMetric}
    For any Hermitian star product $\star$ on $M$ and any fiber metric
    $h$ on $E \longrightarrow M$ and any deformation quantization
    $\bullet$ of $E$ there exists a deformation quantization
    $\boldsymbol{h}$ of $h$ which is unique up to isometry.
\end{theorem}
The first theorem relies heavily on the Serre-Swan theorem and the
fact that algebraic $K_0$-theory is stable under formal deformations
\cite{rosenberg:1996a:pre}. In fact, projections and hence projective
modules can always be deformed in an essentially unique way. The
second statement follows for much more general deformed algebras than
only for star products, see \cite{bursztyn.waldmann:2005b}.
\begin{remark}
    \label{remark:Matter}
    \begin{enumerate}
    \item In case $M$ is symplectic, one has even a rather explicit
        Fedosov-like construction for $\bullet$ and $\boldsymbol{h}$
        in terms of connections, see \cite{waldmann:2002b}.
    \item\label{item:EndEDeformation} It turns out that also
        $\Gamma^\infty(\End(E))$ becomes deformed into an associative
        algebra $(\Gamma^\infty(\End(E))[[\lambda]], \star')$ such
        that $\Gamma^\infty(E)[[\lambda]]$ becomes a Morita
        equivalence bimodule between the two deformed algebras $\star$
        and $\star'$. Together with the deformation $\boldsymbol{h}$
        of $h$ one obtains even a strong Morita equivalence bimodule
        \cite{bursztyn.waldmann:2005b}.
    \item Note also that the results of the two theorems are more than
        just the `analogy' used in the more general framework of
        noncommutative geometry: we have here a precise link between
        the noncommutative geometries and their classical/commutative
        limits via deformation. For general noncommutative geometries
        it is not even clear what a classical/commutative limit is.
    \end{enumerate}
\end{remark}

%
% Deformed principal bundles
%

\section{Deformed principal bundles}
\label{sec:DeformedPrincipalBundles}

This section contains a review of results obtained in
\cite{bordemann.neumaier.waldmann.weiss:2007a:pre} as well as in
\cite{weiss:2006a}.

In all fundamental theories of particle physics the field theories
involve gauge fields. Geometrically, their formulation is based on the
use of a principal bundle $\pr: P \longrightarrow M$ with structure
group $G$, i.e. $P$ is endowed with a (right) action of $G$ which is
proper and free whence the quotient $P \big/ G = M$ is again a smooth
manifold. All the matter fields are then obtained as sections of
associated vector bundles by choosing an appropriate representation of
$G$.

In the noncommutative framework there are several approaches to gauge
theories: for particular structure groups and representations notions
of gauge theories have been developed by Jurco, Schupp, Wess and
coworkers \cite{jurco.et.al.:2001a, jurco.schupp.wess:2001a,
  jurco.schraml.schupp.wess:2000a, jurco.schupp.wess:2000a,
  jurco.schupp:2000a}. Here the focus was mainly on local
considerations and the associated bundles but not on the principal
bundle directly. Conversely, there is a purely algebraic and
intrinsically global formulation of Hopf-Galois extensions where not
only the base manifold $M$ is allowed to be noncommutative but even
the structure group is replaced by a general Hopf algebra, see e.g.
\cite{dabrowski.grosse.hajac:2001a} and references therein for the
relation of Hopf-Galois theory to noncommutative gauge field theories.
However, as we shall see below, in this framework which a priori does
not refer to any sort of deformation, in general only very particular
Poisson structures on $M$ can be used. Finally, in
\cite{waldmann:2002a} a local approach to principal $\mathrm{Gl}(n
\mathbb{C})$ or $\mathrm{U}(n)$ bundles was implicitly used via
deformed transition matrices.

We are now seeking for a definition of a deformation quantization of a
principal bundle $P$ for a generic structure Lie group $G$, arbitrary
$M$ and arbitrary star product $\star$ on $M$ without further
assumptions on $P$. In particular, the formulation should be
intrinsically global.

The idea is to consider the classical algebra homomorphism
\begin{equation}
    \label{eq:prpullback}
    \pr^*: C^\infty(M) \longrightarrow C^\infty(P)
\end{equation}
and try to find a reasonable deformation of $\pr^*$. The first idea
would be to find a star product $\star_P$ on $P$ with a deformation
$\boldsymbol{\pr^*} = \sum_{r=0}^\infty \lambda^r
\boldsymbol{\pr^*}_r$ of $\boldsymbol{\pr^*}_0 = \pr^*$ into an
algebra homomorphism 
\begin{equation}
    \label{eq:prAlgebraMorph}
    \boldsymbol{\pr^*}(f \star g) 
    = \boldsymbol{\pr^*}(f) \star_P \boldsymbol{\pr^*}(g)
\end{equation}
with respect to the two star products $\star$ and $\star_P$. In some
sense this would be the first (but not the only) requirement for a
Hopf-Galois extension. In fact, the first order of
\eqref{eq:prAlgebraMorph} implies that the \emph{classical} projection
map $\pr$ is a Poisson map with respect to the Poisson structures
induced by $\star$ on $M$ and $\star_P$ on $P$. The following example
shows that in general there are obstructions to achieve
\eqref{eq:prAlgebraMorph} already on the classical level:
\begin{example}
    \label{example:HopfFibration}
    Consider the Hopf fibration $\pr: S^3 \longrightarrow S^2$ (which
    is a nontrivial principal $S^1$-bundle over $S^2$) and equip $S^2$
    with the canonical symplectic Poisson structure. Then there exists
    \emph{no} Poisson structure on $S^3$ such that $\pr$ becomes a
    Poisson map. Indeed, if there would be such a Poisson structure
    then necessarily all symplectic leaves would be two-dimensional as
    symplectic leaves are mapped into symplectic leaves and $S^2$ is
    already symplectic. Fixing one symplectic leaf in $S^3$ one checks
    that $\pr$ restricted to this leaf is still surjective and thus
    provides a covering of $S^2$. But $S^2$ is simply connected whence
    the symplectic leaf is itself a $S^2$. This would yield a section
    of the nontrivial principal bundle $\pr: S^3 \longrightarrow S^2$,
    a contradiction.
\end{example}
\begin{remark}
    \label{remark:OtherHopfGalois}
    Note that there are prominent examples of Hopf-Galois extensions
    using quantum spheres, see
    e.g.~\cite{hajac.matthes.szymanski:2003a} and references therein.
    The above example shows that when taking the semi-classical limit
    of these $q$-deformations one obtains Poisson structures on $S^2$
    which are certainly not symplectic.  Note that this was a crucial
    feature in the above example. A further investigation of these
    examples is work in progress.
\end{remark}

The above example shows that the first idea of deforming the
projection map into an algebra homomorphism leads to hard obstructions
in general, even though there are interesting classes of examples
where the obstructions are absent. However, as we are interested in an
approach not making too much assumptions in the beginning, we abandon
this first idea. The next weaker requirement would be to deform
$\pr^*$ not into an algebra homomorphism but only turning
$C^\infty(P)$ into a \emph{bimodule}. This would have the advantage
that there is no Poisson structure on $P$ needed. However, a more
subtle analysis shows that again for the Hopf fibration such a
bimodule structure is impossible if one uses a star product on $S^2$
coming from the symplectic Poisson structure. Thus we are left with a
\emph{module structure}: for later convenience we choose a right
module structure and state the following definition
\cite{bordemann.neumaier.waldmann.weiss:2007a:pre}:
\begin{definition}
    \label{definition:DeformationPFB}
    Let $\pr: P \circlearrowleft G \longrightarrow M$ be a principal
    $G$-bundle over $M$ and $\star$ a star product on $M$. A
    deformation quantization of $P$ is a right $\star$-module
    structure $\bullet$ for $C^\infty(P)[[\lambda]]$ of the form
    \begin{equation}
        \label{eq:BulletPFB}
        F \bullet f = F \pr^*f + \sum_{r=1}^\infty \lambda^r
        \varrho_r(F, f),
    \end{equation}
    where $\varrho_r: C^\infty(P) \times C^\infty(M) \longrightarrow
    C^\infty(P)$ is a bidifferential operator (along $\pr$) for all $r
    \ge 1$, such that in addition one has the $G$-equivariance
    \begin{equation}
        \label{eq:Gequivariance}
        g^* (F \bullet f) = g^*F \bullet f
    \end{equation}
    for all $F \in C^\infty(P)[[\lambda]]$, $f \in
    C^\infty(M)[[\lambda]]$ and $g \in G$.
\end{definition}
Note that as $G$ acts on $P$ from the right, the pull-backs with the
actions of $g \in G$ provide a left action on $C^\infty(P)$ in
\eqref{eq:Gequivariance}. Then this condition means that the
$G$-action commutes with the module multiplications.

Note that the module property $F \bullet (f \star g) = (F \bullet f)
\bullet g$ implies that the constant function $1$ acts as identity.
Indeed, since $1 \star 1 = 1$ the action of $1$ via $\bullet$ is a
projection. However, in zeroth order the map $F \mapsto F \bullet 1$
is just the identity and hence invertible. But the only invertible
projection is the identity map itself. Thus
\begin{equation}
    \label{eq:FbulletEins}
    F \bullet 1 = F
\end{equation}
for all $F \in C^\infty(P)[[\lambda]]$, so the module structure
$\bullet$ is necessarily unital.

Finally, we call two deformation quantizations $\bullet$ and
$\tilde{\bullet}$ \emph{equivalent}, if there exists a $G$-equivariant
equivalence transformation between them, i.e. a formal series of
differential operators $T = \id + \sum_{r=1}^\infty \lambda^r T_r$ on
$C^\infty(P)[[\lambda]]$ such that
\begin{equation}
    \label{eq:EquivalenceG}
    T(F \bullet f) = T(F) \tilde{\bullet} f
    \quad
    \textrm{and}
    \quad
    g^* T = T g^*
\end{equation}
for all $F \in C^\infty(P)[[\lambda]]$, $f \in C^\infty(M)[[\lambda]]$
and $g \in G$.

We shall now discuss the existence and classification of such module
structures. For warming up we consider the situation of a
\emph{trivial} principal fiber bundle:
\begin{example}
    \label{example:TrivialPFB}
    Let $P = M \times G$ be the trivial (and trivialized) principal
    $G$-bundle over $M$ with the obvious projection. For any star
    product $\star$ on $M$ we can now extend $\star$ to $C^\infty(M
    \times G)[[\lambda]]$ by simply acting only on the $M$-coordinates
    in the Cartesian product. Here we use the fact that we can
    canonically extend multidifferential operators on $M$ to $M \times
    G$. Clearly, all algebraic properties are preserved whence in this
    case we even get a star product $\star_P = \star \otimes \mu$ with
    the undeformed multiplication $\mu$ for the $G$-coordinates. In
    particular, $C^\infty(M \times G)[[\lambda]]$ becomes a right
    module with respect to $\star$. So locally there are no
    obstructions even for the strongest requirement
    \eqref{eq:prAlgebraMorph} and hence also for \eqref{eq:BulletPFB}.
\end{example}

The problem of finding $\bullet$ is a global question whence we can
not rely on local considerations directly. The most naive way to
construct a $\bullet$ is an order-by-order construction: In general,
one has to expect obstructions in each order which we shall now
compute explicitly. This is a completely standard approach from the
very first days of algebraic deformation theory
\cite{gerstenhaber:1964a} and will in general only yield the result
that there are possible obstructions: in this case one needs more
refined arguments to ensure existence of deformations whence the
order-by-order argument in general is rather useless. In our situation,
however, it turns out that we are surprisingly lucky.

The following argument applies essentially to arbitrary algebras and
module deformations and should be considered to be folklore.  Suppose
we have already found $\varrho_0 = \pr^*$, $\varrho_1$, \ldots,
$\varrho_k$ such that 
\begin{equation}
    \label{eq:bulletk}
    F \bullet^{(k)} f 
    = F \pr^* f + \sum_{r=1}^k \lambda^r \varrho_r(F, f)
\end{equation}
is a module  structure up to order $\lambda^k$ and each $\varrho_r$
fulfills the $G$-equivariance condition. Then in order to find
$\varrho_{k+1}$ such that $\bullet^{(k+1)} = \bullet^{(k)} +
\lambda^{k+1} \varrho_{k+1}$ is a module structure up to order
$\lambda^{k+1}$ we have to satisfy
\begin{align}
    &\varrho_{k+1}(F, f) \pr^*g 
    - \varrho_{k+1}(F, fg) +
    \varrho_{k+1}(F \pr^*f, g)
    \nonumber \\
    &\qquad=
    \sum_{r=1}^k 
    \left(
        \varrho_r(F, C_{k+1-r}(f, g)) 
        - 
        \varrho_r(\varrho_{k+1-r}(F, f), g)
    \right)
    = R_k (F, f, g),
    \label{eq:varrhokpluseins}
\end{align}
for all $F \in C^\infty(P)[[\lambda]]$ and $f, g \in
C^\infty(M)[[\lambda]]$. Here $C_r$ denotes the $r$-th cochain of the
star product $\star$ as in \eqref{eq:StarProduct}. In order to
interpret this equation we consider the $\varrho_r$ as maps
\begin{equation}
    \label{eq:varrhorCochains}
    \varrho_r: C^\infty(M) \ni f
    \; \mapsto \;
    \varrho_r(\cdot, f) \in \Diffop(P)
\end{equation}
and similarly
\begin{equation}
    \label{eq:Rk}
    R_k: C^\infty(M) \times C^\infty(M) \ni (f, g)
    \; \mapsto \;
    R_k( \cdot, f, g) \in \Diffop(P).
\end{equation}
Viewing $\Diffop(P)$ as $C^\infty(M)$-bimodule via $\pr^*$ in the
usual way, we can now re-interpret \eqref{eq:varrhokpluseins} as
equation between a Hochschild one-cochain $\varrho_{k+1}$ and a
Hochschild two-cochain $R_k$
\begin{equation}
    \label{eq:deltavarrhokpluseins}
    \delta \varrho_{k+1} = R_k
\end{equation}
in the Hochschild (sub-)complex $\HCdiff^\bullet (C^\infty(M),
\Diffop(P))$ consisting of \emph{differential} cochains taking values
in the bimodule $\Diffop(P)$. Here $\delta$ is the usual Hochschild
differential. Using the assumption that the $\varrho_0, \ldots,
\varrho_k$ have been chosen such that $\bullet^{(k)}$ is a module
structure up to order $\lambda^k$ it is a standard argument to show
\begin{equation}
    \label{eq:deltaRkNull}
    \delta R_k = 0.
\end{equation}
Thus the necessary condition for \eqref{eq:deltavarrhokpluseins} is
always fulfilled by construction whence
\eqref{eq:deltavarrhokpluseins} is a cohomological condition: The
equation \eqref{eq:deltavarrhokpluseins} has solutions if and only if
the class of $R_k$ in the second Hochschild cohomology
$\HHdiff^2(C^\infty(M), \Diffop(P))$ is trivial.

In fact, we have also to take care of the $G$-equivariance of
$\varrho_{k+1}$. If all the $\varrho_0$, \ldots, $\varrho_k$ satisfy
the $G$-equivariance then it is easy to see that also $R_k$ has the
$G$-equivariance property. Thus we have to consider yet another
subcomplex of the differential Hochschild complex, namely
\begin{equation}
    \label{eq:HCdiffG}
    \HCdiff^\bullet (C^\infty(M), \Diffop(P)^G)
    \subseteq
    \HCdiff^\bullet (C^\infty(M), \Diffop(P)).
\end{equation}
Thus the obstruction for \eqref{eq:deltavarrhokpluseins} to have a
$G$-equivariant solution is the Hochschild cohomology class
\begin{equation}
    \label{eq:RkinHH}
    [R_k] \in \HHdiff^2 (C^\infty(M), \Diffop(P)^G).
\end{equation}

A completely analogous order-by-order construction shows that also the
obstructions for equivalence of two deformations $\bullet$ and
$\tilde{\bullet}$ can be formulated using the differential Hochschild
complex of $C^\infty(M)$ with values in $\Diffop(P)^G$. Now the
obstruction lies in the first cohomology $\HHdiff^1(C^\infty(M),
\Diffop(P)^G)$.

The following (nontrivial) theorem solves the problem of existence and
uniqueness of deformation quantizations now in a trivial way
\cite{bordemann.neumaier.waldmann.weiss:2007a:pre}:
\begin{theorem}
    \label{theorem:HHNull}
    Let $\pr: P \longrightarrow M$ be a surjective submersion.
    \begin{enumerate}
    \item We have
        \begin{equation}
            \label{eq:HHdiffNull}
            \HHdiff^k(C^\infty(M), \Diffop(P)) 
            = 
            \begin{cases}
                \Diffop_\ver(P) & \textrm{for} \; k=0 \\
                \{0\} & \textrm{for} \; k \ge 1.
            \end{cases}
        \end{equation}
    \item If in addition $\pr: P \circlearrowleft G \longrightarrow M$
        is a principal $G$-bundle then we have
        \begin{equation}
            \label{eq:HHdiffGNull}
            \HHdiff^k(C^\infty(M), \Diffop(P)^G) 
            = 
            \begin{cases}
                \Diffop_\ver(P)^G & \textrm{for} \; k=0 \\
                \{0\} & \textrm{for} \; k \ge 1.
            \end{cases}
        \end{equation}
    \end{enumerate}
\end{theorem}
The main idea is to proceed in three steps: first one shows that one
can localize the problem to a bundle chart. For the local situation
one can use the explicit homotopies from
\cite{bordemann.et.al:2005a:pre} to show that the cohomology is
acyclic. This is the most nontrivial part. By a suitable partition of
unity one can glue things together to end up with the global
statement. For a detailed proof we refer to
\cite{bordemann.neumaier.waldmann.weiss:2007a:pre}.

From this theorem and the previous considerations we obtain
immediately the following result
\cite{bordemann.neumaier.waldmann.weiss:2007a:pre}:
\begin{corollary}
    \label{corollary:DeformPFBOk}
    For every principal $G$-bundle $\pr: P \circlearrowleft G
    \longrightarrow M$ and any star product $\star$ on $M$ there
    exists a deformation quantization $\bullet$ which is unique up to
    equivalence.
\end{corollary}
In particular, the deformation for the trivial bundle as in
Example~\ref{example:TrivialPFB} is the unique one up to equivalence.

\begin{remark}
    \label{remark:DeformationPFB}
    \begin{enumerate}
    \item It should be noted that the use of
        Theorem~\ref{theorem:HHNull} gives existence and uniqueness
        but no explicit construction of deformation quantizations of
        principal bundles. Here the cohomological method is not
        sufficient even though in
        \cite{bordemann.neumaier.waldmann.weiss:2007a:pre} rather
        explicit homotopies were constructed which allow to determine
        further properties of $\bullet$.
    \item In the more particular case of a symplectic Poisson
        structure on $M$, Weiss used in his thesis \cite{weiss:2006a}
        a variant of Fedosov's construction which gives a much more
        geometric and explicit approach: there is a well-motivated
        geometric input, namely a symplectic covariant derivative on
        $M$ as usual for Fedosov's star products and a principal
        connection on $P$. Out of this the module multiplication
        $\bullet$ is constructed by a recursive procedure. The
        dependence of $\bullet$ on the principal connection should be
        interpreted as a global version of the Seiberg-Witten map
        \cite{seiberg.witten:1999a}, now of course in a much more
        general framework for arbitrary principal bundles, see also
        \cite{barnich.brandt.grigoriev:2002a,
          jurco.schraml.schupp.wess:2000a, jurco.et.al.:2001a}.
    \item For the general Poisson case a more geometric construction
        is still missing. However, it seems to be very promising to
        combine global formality theorems like the one in
        \cite{dolgushev:2005a} or the approach in
        \cite{cattaneo.felder.tomassini:2002b} with the construction
        \cite{weiss:2006a}. These possibilities will be investigated
        in future works.
    \end{enumerate}
\end{remark}

%
% The commutant and associated bundles
%

\section{The commutant and associated bundles}
\label{sec:CommutantAssociated}

Theorem~\ref{theorem:HHNull} gives in addition to the existence and
uniqueness of deformation quantizations of $P$ also a description of
the \emph{differential commutant} of the right multiplications by
functions on $M$ via $\bullet$: we are interested in those formal
series $D = \sum_{r=0}^\infty \lambda^r D_r \in \Diffop(P)[[\lambda]]$
of differential operators with the property
\begin{equation}
    \label{eq:Commutes}
    D(F \bullet f) = D(F) \bullet f
\end{equation}
for all $F \in C^\infty(P)[[\lambda]]$ and $f \in
C^\infty(M)[[\lambda]]$. In particular, if $D_0 = \id$ then
\eqref{eq:Commutes} gives a \emph{self-equivalence}. Clearly, the
differential commutant
\begin{equation}
    \label{eq:Kommutante}
    \boldsymbol{\mathcal{K}} = 
    \left\{
        D \in \Diffop(P)[[\lambda]] 
        \; \big| \; 
        D \; \textrm{satisfies} \; \eqref{eq:Commutes}
    \right\}
    \subseteq \Diffop(P)[[\lambda]]
\end{equation}
is a subalgebra of $\Diffop(P)[[\lambda]]$ over
$\mathbb{C}[[\lambda]]$.

Note that there are other operators on $C^\infty(P)[[\lambda]]$ which
commute with all right multiplications, namely the highly non-local
pull-backs $g^*$ with $g \in G$. This was just part of the
Definition~\ref{definition:DeformationPFB} of a deformation
quantization of a principal bundle. However, in this section we shall
concentrate on the differential operators with \eqref{eq:Commutes}
only.

Before describing the commutant it is illustrative to consider the
classical situation. Here the commutant is simply given by the
vertical differential operators
\begin{equation}
    \label{eq:ClassicalKommutant}
    \Diffop_\ver(P) = 
    \left\{
        D \in \Diffop(P)
        \; \big| \; 
        D(F \pr^*f) = D(F) \pr^*f
    \right\}
\end{equation}
by the very definition of vertical differential operators.
Alternatively, the commutant is the zeroth Hochschild cohomology. More
interesting is now the next statement which gives a quantization of
the classical commutant, see
\cite{bordemann.neumaier.waldmann.weiss:2007a:pre}.
\begin{theorem}
    \label{theorem:Commutant}
    There exists a $\mathbb{C}[[\lambda]]$-linear bijection
    \begin{equation}
        \label{eq:varrhoprime}
        \varrho': \Diffop_\ver(P)[[\lambda]] 
        \longrightarrow
        \boldsymbol{\mathcal{K}} 
        \subseteq \Diffop(P)[[\lambda]]
    \end{equation}
    of the form
    \begin{equation}
        \label{eq:varrhoprimeExplicit}
        \varrho' = \id + \sum_{r=1}^\infty \lambda^r \varrho'_r
    \end{equation}
    which is $G$-equivariant, i.e.
    \begin{equation}
        \label{eq:Gequivariantvarrhoprime}
        g^* \varrho' = \varrho' g^*
    \end{equation}
    for all $g \in G$. The choice of such a $\varrho'$ induces an
    associative deformation $\star'$ of $\Diffop_\ver(P)[[\lambda]]$
    which is uniquely determined by $\star$ up to
    equivalence. Finally, $\varrho'$ induces a left
    $(\Diffop_\ver(P)[[\lambda]], \star')$-module structure $\bullet'$
    on $C^\infty(P)[[\lambda]]$ via
    \begin{equation}
        \label{eq:bulletprimeDef}
        D \bullet' F = \varrho'(D) F.
    \end{equation}
\end{theorem}
The proof relies on an adapted symbol calculus for the differential
operators $\Diffop(P)$: using an appropriate $G$-invariant covariant
derivative $\nabla^P$ on $P$ which preserves the vertical distribution
and a principal connection on $P$ one can induce a $G$-equivariant
splitting of the differential operators $\Diffop(P)$ into the vertical
differential operators and those differential operators which
differentiate at least once in horizontal directions. Note that this
complementary subspace has no intrinsic meaning but depends on the
choice of $\nabla^P$ and the principal connection. A recursive
construction gives the corrections terms $\varrho'_r(D)$ for a given
$D \in \Diffop_\ver(P)$, heavily using the fact that the \emph{first}
Hochschild cohomology $\HHdiff^1(C^\infty(M), \Diffop(P))$ vanishes.
Since the commutant itself is an associative algebra the remaining
statements follow.
\begin{corollary}
    \label{corollary:starprimeInvariant}
    For the above choice of $\varrho'$ the resulting deformation
    $\star'$ as well as the module structure are $G$-invariant, i.e. we
    have
    \begin{equation}
        \label{eq:Ginvariance}
        g^*(D \star' \tilde{D}) = g^*D \star' g^*\tilde{D}
        \quad
        \textrm{and}
        \quad
        g^*(D \bullet' F) = g^*D \bullet' g^*F
    \end{equation}
    for all $D, \tilde{D} \in \Diffop_\ver(P)[[\lambda]]$ and $F \in
    C^\infty(P)[[\lambda]]$.
\end{corollary}
This follows immediately from the $G$-equivariance of $\bullet$ and
the $G$-equivariance of $\varrho'$.
\begin{remark}
    \label{remark:Bicommutant}
    A simple induction shows that the commutant of
    $(\Diffop_\ver(P)[[\lambda]], \star')$ inside all differential
    operators $\Diffop(P)[[\lambda]]$ is again
    $(C^\infty(M)[[\lambda]], \star)$, where both algebras act by
    $\bullet'$ and $\bullet$, respectively. This way
    $C^\infty(P)[[\lambda]]$ becomes a $(\star', \star)$-bimodule such
    that the two algebras acting from left and right are mutual
    commutants inside all differential operators. Though this
    resembles already much of a Morita context, it is easy to see that
    $C^\infty(P)[[\lambda]]$ is \emph{not} a Morita equivalence
    bimodule, e.g it is not finitely generated and projective.
    However, as we shall see later, there is still a close relation to
    Morita theory to be expected.
\end{remark}
\begin{remark}
    \label{remark:BimoduleDeformation}
    Note that classically $\pr^*: C^\infty(M) \longrightarrow
    \Diffop(P)$ is an algebra homomorphism, too. Thus the questions
    raised at the beginning of
    Section~\ref{sec:DeformedPrincipalBundles} can now be rephrased as
    follows: for a \emph{bimodule deformation} of $C^\infty(P)$ into a
    bimodule over $C^\infty(M)[[\lambda]]$ equipped with possibly two
    different star products for the left and right action, one has to
    deform $\pr^*$ into a map
    \begin{equation}
        \label{eq:deformBimodule}
        \boldsymbol{\pr^*}:
        C^\infty(M)[[\lambda]] \longrightarrow
        (\Diffop_\ver(P)[[\lambda]], \star')
    \end{equation}
    such that the image is a subalgebra. In this case, we can induce a
    new product $\star'_M$ also for $C^\infty(M)[[\lambda]]$ making
    $C^\infty(P)[[\lambda]]$ a bimodule for the two, possibly
    different, star product algebras $(C^\infty(M)[[\lambda]],
    \star_M')$ from the left and $(C^\infty(M)[[\lambda]], \star)$ from
    the right. Note that this is the only way to achieve it since
    $\star'$ is uniquely determined by $\star$. Thus it is clear that
    we have to expect obstructions in the general case as there might
    be no subalgebra of $(\Diffop_\ver(P)[[\lambda]], \star')$ which
    is in bijection to $C^\infty(M)[[\lambda]]$. Even if this might be
    the case, the resulting product $\star'_M$ might be inequivalent
    to $\star$.  Note however, that we have now a very precise
    framework for the question whether $\pr^*$ can be deformed into a
    bimodule structure.
\end{remark}
\begin{remark}
    \label{remark:DeftoDef}
    As a last remark we note that changing $\star$ to an equivalent
    $\tilde{\star}$ via an equivalence transformation $\Phi$ yields a
    corresponding right module structure $\tilde{\bullet}$ by
    \begin{equation}
        \label{eq:equivalentstarbullet}
        F \tilde{\bullet} f =  F \bullet \Phi(f),
    \end{equation}
    which is still unique up to equivalence by
    Theorem~\ref{theorem:HHNull}. It follows that the commutants are
    \emph{equal} (for this particular choice of $\tilde{\bullet}$)
    whence the induced deformations $\star'$ and $\tilde{\star}'$
    coincide. An equivalent choice of $\tilde{\bullet}$ would result
    in an equivalent $\tilde{\star}'$. This shows that we obtain a
    well-defined map
    \begin{equation}
        \label{eq:DeftoDef}
        \Def(C^\infty(M)) 
        \longrightarrow 
        \Def(\Diffop_\ver(P))
    \end{equation}
    for the sets of equivalence classes of associative deformations.
    In fact, the resulting deformations $\star'$ are even
    $G$-invariant, whence the above map takes values in the smaller
    class of $G$-invariant deformations $\Def_G(\Diffop_\ver(P))$.
\end{remark}

To make contact with the deformed vector bundles from
Section~\ref{sec:MatterFieldsVectorBundles} we consider now the
association process. Recall that on the classical level one starts
with a (continuous) representation $\pi$ of $G$ on a
finite-dimensional vector space $V$. Then the associated vector bundle
is
\begin{equation}
    \label{eq:AssoVec}
    E = P \times_G V \longrightarrow M,
\end{equation}
where the fibered product is defined via the equivalence relation $(p
\cdot g, v) \sim (p, \pi(g) v)$ as usual. As the action of $G$ on $P$
is proper and free, $E$ is a smooth manifold again and, in fact, a
vector bundle over $M$ with typical fiber $V$. Rather tautologically,
any vector bundle is obtained like this by association from its own
frame bundle. For the sections of $E$ one has the canonical
identifications
\begin{equation}
    \label{eq:SectionE}
    \Gamma^\infty(E) \cong C^\infty(P, V)^G
\end{equation}
as right $C^\infty(M)$-modules, where the $G$-action of $C^\infty(P,
V)$ is the obvious one.

After this preparation it is clear how to proceed in the deformed
case. From the $G$-equivariance of $\bullet$ we see that
\begin{equation}
    \label{eq:GammaE}
    \Gamma^\infty(E)[[\lambda]] \cong C^\infty(P, V)^G[[\lambda]]
    \subseteq C^\infty(P, V)[[\lambda]]
\end{equation}
is a $\star$-submodule with respect to the restricted module
multiplication $\bullet$. It induces a right $\star$-module structure
for $\Gamma^\infty(E)[[\lambda]]$ which we still denote by $\bullet$.
This way we recover the deformed vector bundle as in
Section~\ref{sec:MatterFieldsVectorBundles}.

Moreover, we see that the $\End(V)$-valued differential operators
$\Diffop(P) \otimes \End(V)$ canonically act on $C^\infty(P, V)$
whence $\left((\Diffop_\ver(P) \otimes \End(V))[[\lambda]],
    \star'\right)$ acts via $\bullet'$ on $C^\infty(P, V)[[\lambda]]$
in such a way that the action commutes with the
$\bullet$-multiplications from the right. By the $G$-invariance of
$\star'$ we see that the invariant elements
$\left(\Diffop_\ver(P)\otimes\End(V)\right)^G[[\lambda]]$ form a
$\star'$-subalgebra which preserves (via $\bullet'$) the
$\bullet$-submodule $C^\infty(P, V)^G[[\lambda]]$. Thus we obtain an
algebra homomorphism
\begin{equation}
    \label{eq:MotherOfAll}
    \left(
        (\Diffop_\ver(P) \otimes \End(V))^G[[\lambda]], \star'
    \right)
    \longrightarrow
    \left(
        \Gamma^\infty(\End(E))[[\lambda]], \star'
    \right)
\end{equation}
where $\star'$ on the left hand side is the deformation from
Remark~\ref{remark:Matter}, part~\ref{item:EndEDeformation}.

We conclude this section with some remarks and open questions:
\begin{remark}
    \label{remark:OpenQuestions}
    \begin{enumerate}
    \item The universal enveloping algebra valued gauge fields of
        \cite{jurco.et.al.:2001a, jurco.schraml.schupp.wess:2000a} can
        now easily be understood. For two vertical \emph{vector
          fields} $\xi, \eta \in \Diffop_\ver(P)$ we have an action on
        $C^\infty(P)[[\lambda]]$ via $\bullet'$-left multiplication.
        In zeroth order this is just the usual Lie derivative
        $\Lie_\xi$. Now the module structure says that
        \begin{equation}
            \label{eq:Ugvalued}
            \xi \bullet' (\eta \bullet' F) 
            - \eta \bullet' (\xi \bullet' F)
            =
            ([\xi, \eta]_{\star'}) \bullet' F
        \end{equation}
        for all $F \in C^\infty(P)[[\lambda]]$. Here $[\xi,
        \eta]_{\star'} = \xi \star' \eta - \eta \star' \xi \in
        \Diffop_\ver(P)[[\lambda]]$ is the $\star'$-commutator. In
        general, this commutator is a formal series of vertical
        differential operators but not necessarily a vector field any
        more. Note that \eqref{eq:Ugvalued} holds already on the level
        of the principal bundle.
    \item For noncommutative gauge field theories we still need a good
        notion of gauge fields, i.e. connection one-forms, and their
        curvatures within our global approach. Though there are
        several suggestions from e.g. \cite{jurco.schupp.wess:2000a} a
        conceptually clear picture seems still to be missing.
    \item In a future project we plan to investigate the precise
        relationship between $(\Diffop_\ver(P)[[\lambda]], \star')$
        and the Morita theory of star products
        \cite{bursztyn.waldmann:2001a, bursztyn.waldmann:2002a,
          bursztyn.waldmann:2000b}. Here \eqref{eq:MotherOfAll}
        already suggests that one can re-construct all algebras Morita
        equivalent to $(C^\infty(M)[[\lambda]], \star)$ out of
        $\star'$.
    \end{enumerate}
\end{remark}

%
% Acknowledgment
%

\subsection*{Acknowledgment}

It is a pleasure for me to thank the organizers Bertfried Fauser,
Jürgen Tolksdorf and Eberhard Zeidler for their invitation to the very
stimulating conference ``Recent Developments in Quantum Field
Theory''. Moreover, I would like to thank Rainer Matthes for valuable
discussions on Hopf-Galois extensions and Stefan Weiß for many
comments on the first draft of the manuscript.

%
% references
%

%\begin{footnotesize}
%    \renewcommand{\arraystretch}{0.5} 
%    \bibliographystyle{ewde}
%    \bibliography{dqarticle,dqbook,dqprocentry,dqproceeding,preprints,misc,dqthesis,notes}
%\end{footnotesize}

\end{document}